\documentclass[12pt,oneside,reqno]{amsart}
\usepackage{}
\usepackage{amssymb}
\usepackage{bbm}
\usepackage{cases}
\usepackage{amsmath}
\usepackage{graphicx}
\usepackage{mathrsfs}
\usepackage{stmaryrd}
\usepackage{amsfonts}
\usepackage{enumerate,amsmath,amssymb,amsthm}

\numberwithin{equation}{section}

\newcommand{\be}{\begin{eqnarray}}
\newcommand{\ee}{\end{eqnarray}}
\newcommand{\ce}{\begin{eqnarray*}}
\newcommand{\de}{\end{eqnarray*}}
\newtheorem{theorem}{Theorem}[section]
\newtheorem{lemma}[theorem]{Lemma}
\newtheorem{remark}[theorem]{Remark}
\newtheorem{definition}[theorem]{Definition}
\newtheorem{proposition}[theorem]{Proposition}
\newtheorem{Examples}[theorem]{Example}
\newtheorem{corollary}[theorem]{Corollary}

\def\e{{\mathrm{e}}}
\def\eps{\varepsilon}

\def\p{\partial}

\def\[{{\Big[}}
\def\]{{\Big]}}
\def\<{{\langle}}
\def\>{{\rangle}}
\def\({{\Big(}}
\def\){{\Big)}}

\def\bx{{\mathbf{x}}}

\def\dif{{\mathord{{\rm d}}}}

\def\no{\nonumber}
\def\={&\!\!=\!\!&}
\def\bt{\begin{theorem}}
\def\et{\end{theorem}}
\def\bl{\begin{lemma}}
\def\el{\end{lemma}}
\def\br{\begin{remark}}
\def\er{\end{remark}}

\def\bd{\begin{definition}}
\def\ed{\end{definition}}
\def\bp{\begin{proposition}}
\def\ep{\end{proposition}}
\def\bc{\begin{corollary}}
\def\ec{\end{corollary}}
\def\bx{\begin{Examples}}
\def\ex{\end{Examples}}

\def\cJ{{\mathcal J}}

\def\cM{{\mathcal M}}

\def\cO{{\mathcal O}}

\def\cQ{{\mathcal Q}}

\def\cT{{\mathcal T}}

\def\mD{{\mathbb D}}
\def\mE{{\mathbb E}}

\def\mH{{\mathbb H}}
\def\mI{{\mathbb I}}

\def\mL{{\mathbb L}}

\def\mN{{\mathbb N}}

\def\mP{{\mathbb P}}
\def\mQ{{\mathbb Q}}
\def\mR{{\mathbb R}}

\def\sB{{\mathscr B}}

\def\sF{{\mathscr F}}

\def\sL{{\mathscr L}}

\def\geq{\geqslant}
\def\leq{\leqslant}

\allowdisplaybreaks

\begin{document}
\title{Well-posedness for SDEs driven by different type of noises}

\date{}

\author{{Yueling Li},\,\,{Longjie Xie}
and  {Yingchao Xie}
}

\address{Yueling Li: School of Mathematics and Statistics, Jiangsu Normal University,
Xuzhou, Jiangsu 221000, P.R.China\\
Email: lylmath@jsnu.edu.cn
}
\address{Longjie Xie:
School of Mathematics and Statistics, Jiangsu Normal University,
Xuzhou, Jiangsu 221000, P.R.China\\
Email: xlj.98@whu.edu.cn
 }
\address{Yingchao Xie: School of Mathematics and Statistics, Jiangsu Normal University,
Xuzhou, Jiangsu 221000, P.R.China\\
Email: ycxie@jsnu.edu.cn
}

\thanks{This work is supported by the National Natural Science Foundation of China (11271169,11401265) and the Project
Funded by the Priority Academic Program Development of Jiangsu Higher Education Institutions.}

\begin{abstract}
We show the existence and uniqueness of strong solutions for stochastic differential equation driven by partial $\alpha$-stable noise and partial
Brownian noise with singular coefficients. The proof is based on the regularity of degenerate mixed type Kolmogorov equation.

\bigskip

\end{abstract}

\maketitle \rm

\section{Introduction and main results}

In this paper, we consider the following stochastic differential equation (SDE) in $\mR^{d_1+d_2}$:
\begin{equation} \label{sde1}
\left\{ \begin{aligned}
&\dif X_t=F(t,X_t,Y_t)\dif t+\dif L_t,\quad X_0=x\in\mR^{d_1},\\
&\dif Y_t=G(t,X_t,Y_t)\dif t+\dif W_t,\quad Y_0=y\in\mR^{d_2},
\end{aligned} \right.
\end{equation}
where $d_1,d_2\geq 1$, $F: \mR^+\times\mR^{d_1+d_2}\rightarrow\mR^{d_1}$ and $G: \mR^+\times\mR^{d_1+d_2}\rightarrow\mR^{d_2}$ are measurable functions, $L_t$ is a $d_1$-dimensional rotationally symmetric
$\alpha$-stable process with $\alpha>1$ and $W_t$ is a $d_2$-dimensional standard Brownian motion both defined on some filtered probability space
$(\Omega,\sF,(\sF_t)_{t\geq 0},\mP)$. SDE (\ref{sde1}) can be seen as SDE in higher dimensional with partial degenerate noises.
In fact, if we let $z:=(x,y)\in\mR^{d_1+d_2}$ and
\begin{equation*}
Z_t:=(X_t,Y_t),
\end{equation*}
define the vector field
\begin{equation*}
B(t,z)=\left(
  \begin{array}{c}
    F(t,z) \\
    G(t,z) \\
  \end{array}
\right),
\end{equation*}
and the matrix
\begin{equation*}
Q=\left(
  \begin{array}{c}
    \mI_{d_1} \\
    0 \\
  \end{array}
\right),
\quad R=\left(
  \begin{array}{c}
    0 \\
    \mI_{d_2} \\
  \end{array}
\right),
\end{equation*}
where $\mI_d$ denotes the $d\times d$ identity matrix. Then, SDE (\ref{sde1}) can be written as
\begin{align}
\dif Z_t=B(t,Z_t)\dif t+Q\dif L_t+R\dif W_t,\quad Z_0=z\in\mR^{d_1+d_2}.    \label{sde2}
\end{align}
SDEs with general L\'evy noise have been intensively studied due to their wide range of applications, see \cite{A-B-W} and references therein.
Note that in the case $Q=R\equiv0$, SDE (\ref{sde2}) is just the ordinary differential equation (ODE)
\begin{align}
Z'(t)=B\big(t,Z(t)\big),\quad Z(0)=z_0\in\mR^{d_1+d_2}.    \label{ode}
\end{align}
Thus, SDE (\ref{sde1}) can also be regarded as stochastic perturbations to (\ref{ode}) by partial L\'evy type noise and partial Brownian noise.

Our aim is to show that noises have regularization effects on the deterministic system. To be more precise, we provide existence and uniqueness
results for SDE (\ref{sde1}) under conditions on $F$ and $G$ which are forbidden for the deterministic ODE (\ref{ode}).
Moreover, our result shows how much effects can noises of different type play in stabling systems.

\vspace{2mm}
The following is the main result of our paper.

\bt\label{main1}
Assume that for $T>0$ and $p,q\in(1,\infty)$ satisfying
\begin{align}
\frac{d_1}{\alpha p}+\frac{d_2}{2p}+\frac{1}{q}<\frac{1}{2},  \label{index}
\end{align}
and $\beta\in(0,1)$ with
\begin{align}
\beta>1-\frac{\alpha}{2}, \label{index2}
\end{align}
we have
$$
F, G\in L^q\Big([0,T];L^p\big(\mR^{d_2}_y;H^\beta_p(\mR^{d_1}_x)\big)\Big).
$$
Then, SDE (\ref{sde1}) admits a unique strong solution for every starting point $(x,y)\in\mR^{d_1+d_2}$.
\et

\vspace{1mm}
Notice that under the above conditions on the drift coefficients $F$ and $G$, the ODE (\ref{ode}) is far from being well-posed,
which in turn implies that noises can help to stably the deterministic systems (see \cite{F-G-P}).
According to the general scheme, to take care of SDEs with singular coefficients, one needs to study the regularity properties
of the solutions to the associated Kolmogorov equation, which in our case is a equation of the following partial degenerate type:
\begin{align}
\p_tu(t,x,y)&-\tfrac{1}{2}\Delta_yu(t,x,y)-\Delta^{\frac{\alpha}{2}}_xu(t,x,y)\no\\
&-F(t,x,y)\cdot\nabla_xu(t,x,y)-G(t,x,y)\cdot\nabla_y u(t,x,y)=f(t,x,y).\label{aa}
\end{align}
We shall prove the optimal regularity result for (\ref{aa}) in Bessel spaces by using interpolation techniques. This has independent interests. We also note that the drift
is only $L^p$ in $y$-variable in which the Brownian noise acts, and has certain fractional Sobolev regularity in $x$-variable
in which the L\'evy noise acts. In fact, the main point is to take care of the balance between the regularities of the mixing noises of continuous Brownian motion and the pure jump L\'evy noise, and seek the minimal integrability index on the coefficients, which improves certain results even in the non-degenerate noise cases, as we shall see in Remark \ref{rem} below.

\vspace{3mm}
By localization technique, we can prove the following uniqueness of strong solution under some local conditions on the coefficients.

\bt\label{main11}
Suppose that for $T>0$, $p,q\in(1,\infty]$ and $\beta\in(0,1)$ satisfying (\ref{index}) and (\ref{index2}) respectively, we have
$$
F, G\in L^q\Big([0,T];L^p_{loc}\big(\mR^{d_2}_y;H^{\beta}_{p,loc}(\mR^{d_1}_x)\big)\Big).
$$
Then, for every starting point $(x,y)\in\mR^{d_1+d_2}$, SDE (\ref{sde1}) has a unique strong solution $(X_t(x),Y_t(y))$ up to the explosion time $\zeta(x,y)$.
\et

The advantage of Theorem \ref{main11} is that the local conditions imposed on the coefficients alow functions with certain growth at infinity,
and we would obtain a unique global strong solution once we can show the unique local solution which does not explode. As an application,
we consider the following stochastic perturbations of Kinetic equation:
\begin{equation} \label{sde3}
\left\{ \begin{aligned}
&\dif X_t=V_t\dif t+\eps\dif L_t,\qquad\quad\,\,\, X_0=x\in\mR^{d_1},\\
&\dif V_t=G(X_t,V_t)\dif t+\dif W_t,\quad V_0=v\in\mR^{d_2},
\end{aligned} \right.
\end{equation}
where $\eps>0$ is a fixed constant. The following result follows directly from Theorem \ref{main11}.

\bc\label{main2}
Assume that for $\beta\in(0,1)$ satisfying (\ref{index2}) and $p\in(1,\infty]$ satisfying
\begin{align*}
\frac{d_1}{\alpha p}+\frac{d_2}{2p}<1-\frac{1}{\alpha} ,
\end{align*}
we have
$$G\in L^p\big(\mR^{d_2}_y;H^{\beta}_p(\mR^{d_1}_x)\big).
$$
Then, for every starting point $(x,v)\in\mR^{d_1+d_2}$,
SDE (\ref{sde3}) has a unique local strong solution.
\ec

Another interesting example of our results is the following SDE driven by $\alpha$-stable noise and absolutely continuous Gaussian-type process:
\begin{align}
\dif X_t=\sigma_1(t,W_t)\dif t+b_1(t,X_t)\dif t+\dif L_t,\quad X_0=x\in\mR^{d_1}.    \label{sdea}
\end{align}
Taking $G\equiv0$ and $F(t,x,y)=b_1(t,x)+\sigma_1(t,y)$ in SDE (\ref{sde1}),
we have:

\bc\label{main3}
Let $d_1,d_2\geq 1$. Assume that for $T>0$,
$$
\sigma_1\in L^q\big([0,T];L^p(\mR^{d_2})\big),\quad b_1\in L^q\big([0,T];H^\beta_p(\mR^{d_1})\big),
$$
where $p,q$ and $\beta$ satisfy (\ref{index}) and (\ref{index2}), respectively.
Then, SDE (\ref{sdea}) has a unique strong solution for every starting point $x\in\mR^{d_1}$.
\ec

Similarly, we can exchange the position of $W_t$ and $L_t$ in (\ref{sdea}) to get that:

\bc\label{main4}
Let $d_1,d_2\geq 1$. Assume that for $T>0$,
$$
\sigma_2\in L^q\big([0,T];H^\beta_p(\mR^{d_1})\big),\quad b_2\in L^q\big([0,T];L^p(\mR^{d_2})\big)
$$
with $p,q$ and $\beta$ satisfying (\ref{index}) and (\ref{index2}), respectively.
Then, SDE
\begin{align}
\dif X_t=\sigma_2(t,L_t)\dif t+b_2(t,X_t)\dif t+\dif W_t,\quad X_0=x\in\mR^{d_2}  \label{sdeb}
\end{align}
admits a unique strong solution for every $x\in\mR^{d_2}$.
\ec

\vspace{1mm}
The regularization effects of noises to the deterministic system have caused much attentions in the past decade.
Consider the following SDE driven by L\'evy noise:
\begin{align}
\dif X_t=b(t,X_t)\dif t+\dif L_t,\quad X_0=x\in\mR^d,\label{sde0}
\end{align}
In the case that $(L_t)_{t\geq 0}$ is a Brownian motion, a remarkable result of Krylov and R\"ockner \cite{Kr-Ro}
shows that SDE (\ref{sde0}) has a unique strong solution if $b\in L^q(\mR_+;L^p(\mR^d))$
with
$$
\frac{d}{p}+\frac{2}{q}<1.
$$
Latter, Zhang \cite{Zh3} extend this result to multiplicative noise under some non-degenerate and Sobolev conditions on the diffusion coefficient.
The case that $(L_t)_{t\geq 0}$ is a pure jump L\'evy process has more difficulties. In fact, when $d=1$ and $(L_t)_{t\geq 0}$ is a symmetric
$\alpha$-stable process with $\alpha<1$, Tanaka, Tsuchiya and Watanabe \cite{Ta-Ts-Wa} showed that even if $b$ is bounded and $\beta$-H\"older
continuous with $\alpha+\beta<1$, SDE (\ref{sde0}) may not has pathwise uniqueness strong solution. On the other hand, when $\alpha\geq 1$ and $b$
is time independent with
\begin{align*}
b\in C_b^\beta(\mR^d),\quad\beta>1-\frac{\alpha}{2},
\end{align*}
it was proved by Priola \cite{Pri} that there exists a unique strong solution $X_t(x)$ for SDE (\ref{sde0}) for each $x\in\mR^d$.
Recently, Zhang \cite{Zh00} obtained the pathwise uniqueness for SDE (\ref{sde0}) when $\alpha>1$ and
$$
b\in L^\infty(\mR^d)\cap H^\beta_p(\mR^d),\quad\beta>1-\frac{\alpha}{2}\quad\text{and}\quad p>2d/\alpha.
$$
In \cite{XZ2}, the authors studied the pathwise uniqueness of singular SDEs driven by general multiplicative L\'evy noise. See also \cite{Ch-So-Zh,Fa-Lu-Th,Fe-Fl-2,Fe-Fl-3,M-N-P-Z,Pri2,W,XZ} for related results and applications. Note that all the works mentioned above are restricted to the non-degenerate noise cases.

In the degenerate noise case, the authors of \cite{F-F-P-V} studied the stochastic Kinetic equation (\ref{sde3}) with $\eps=0$,
and obtained the uniqueness of strong solutions under the condition that
$$
G\in L^p\big(\mR^d;H^\beta_p(\mR^d)\big)\quad\text{with}\quad \beta\in(2/3,1),\,\,p>6d.
$$
See also \cite{W-Z2,W-Z} for H\"older-Dini drift and references therein.

\vspace{2mm}
Compared with the works mentioned above, we would like to make the following comments on our results.

\br\label{rem}
1) Formally, if $d_1=0$ in (\ref{index}), then this means that we require
$$
\frac{d_2}{p}+\frac{2}{q}<1,
$$
which is the same condition on the index as in \cite{Kr-Ro}. On the other hand, if $d_2=0$, then (\ref{index}) becomes that
$$
\frac{d_1}{p}+\frac{\alpha}{q}<\frac{\alpha}{2}.
$$
Compared with \cite{Pri,Zh00}, we drop the boundness condition on drift coefficients in the $x$-variable even in the case that $q=\infty$.
This improvement is certainly un-straightforward, and will be obtained by following the argument in \cite{XZ2} and making use of Sobolev embedding theorem.

2) Note that since $\alpha>1$, it holds that $1-\alpha/2<2/3$. Thus, the condition in Corollary \ref{main2} is weaker than that in \cite{F-F-P-V}. Moreover, the results on SDEs of the form (\ref{sdea}) and (\ref{sdeb}) seem to be new. The main point is that, compared with \cite{F-F-P-V,W-Z2,W-Z}, we can use the regularization effect of L\'evy noise in Corollary \ref{main2} and Corollary \ref{main3} as well as the regularization effect of Brownian motion in Corollary \ref{main4}, thus obtaining weaker conditions on the drift coefficients.
\er

\vspace{2mm}
This paper proceeds as follows. In Section 2, we prepare some notations of spaces and preliminaries inequalities. Section 3 is devoted to the study of
Kolmogorov type equation in Bessel space. Finally, we prove the main results Theorem \ref{main1} in Section 4. Throughout this paper,
we use the convention that $C$ with or without subscripts will denote a positive constant, whose value may change in different places, and
whose dependence on the parameters can be traced from the calculations.

\section{Preliminaries}

First, we introduce some notations. Let $d,p\geq1$, the norm in $L^p(\mR^d)$ will be denoted by $\|\cdot\|_p$.
For $0<r\leq2$, the Bessel potential space $H^{r}_p:=H^{r}_p(\mR^d)$ is given by
$$
H^{r}_p(\mR^d):=\left\{f\in L^p(\mR^d): \Delta^{\frac{r}{2}}f\in L^p(\mR^d)\right\}
$$
with norm
$$
\|f\|_{r,p}:=\|f\|_p+\left\|\Delta^{\frac{r}{2}}f\right\|_p,
$$
where for $0<r<2$, $\Delta^{\frac{r}{2}}$ is the fractional Laplacian operator defined by
$$
\Delta^{\frac{r}{2}}f:=\big(|\cdot|^{r}\hat f\check{\big)},\quad \forall f\in C_0^\infty(\mR^d),
$$
and $\hat f$ (resp. $\check{f}$) denotes the Fourier transform (resp. the Fourier inverse transform) of function $f$.
Note that for $n=1, 2$, $H^{n}_p$ is just the usual Sobolev space with equivalent norm (\cite[p. 135, Theorem 3]{St})
$$
\|f\|_{n,p}=\|f\|_p+\|\nabla^nf\|_p,
$$
here and below, $\nabla$ denotes the weak derivative of $f$. While for $0<r\neq$ integer, the fractional Sobolev space
$W^{r,p}:=W^{r,p}(\mR^d)$ is defined to be the space of functions with
\begin{align*}
\|f\|_{W^{r,p}}:=\|f\|_p+\sum_{k=0}^{[r]}\Bigg(\int_{\mR^d}\!\!\int_{\mR^d}\frac{|\nabla^kf(x)-\nabla^kf(y)|^p}{|x-y|^{d+(r-[r])p}}\dif x\dif y\Bigg)^{1/p}<\infty,
\end{align*}
where $[r]$ denotes the integer part of $r$. In this case, the relation between $H^{r}_p$ and $W^{r,p}$ is that
(cf. \cite[p. 190]{Tri}): for $r>0$, $\eps\in(0,r)$ and $p>1$,
\begin{align}
H^{r+\eps}_p\hookrightarrow  W^{r,p}\hookrightarrow H^{r-\eps}_p,\label{re}
\end{align}
where $A\hookrightarrow B$ denotes that the Banach space $A$ is continuously embedded into the Banach space $B$.
The celebrated Sobolev's embedding theorem (see \cite{Tri}) tells us that:
\begin{itemize}
\item if $r-\frac{d}{p}<0$, then for any $p\leq q\leq pd/(d-rp)$, we have
\begin{align}
H^{r}_p\hookrightarrow L^q(\mR^d).\label{emb2}
\end{align}

\item if $r-\frac{d}{p}>0$ is not an integer, then
\begin{align}
H^{r}_p\hookrightarrow C^{r-\frac{d}{p}}_b(\mR^d),\label{emb}
\end{align}
where for some $r>0$, $C^{r}_b(\mR^d)$ is the usual H\"older space with norm
$$
\|f\|_{C^{r}_b}:=\sum_{i=1}^{[r]}\|\nabla^i f(x)\|_\infty+\big[\nabla^{[r]}f\big]_{r-[r]},
$$
here, for a function $f$ on $\mR^d$ and $\vartheta\in(0,1)$,
\begin{align*}
[f]_{\vartheta}:=\sup_{x\neq y}\frac{|f(x)-f(y)|}{|x-y|^{\vartheta}}.
\end{align*}
\end{itemize}

Given a locally integrable function $f$ on $\mR^d$, the Hardy-Littlewood maximal function of $f$ is defined by
$$
\cM f(x):=\sup_{0<r<\infty}\frac{1}{|B_r|}\int_{B_r}|f(x+y)|\dif y,
$$
where $|B_r|$ denotes the Lebesgue measure of ball $B_r$. The following well known result can be found in \cite[P.5, Theorem 1]{St} and \cite{Zh3}.

\bl
(i) For $p\in(1,\infty]$ and all $f\in L^p(\mR^d)$, there exists a constant $C_{d,p}>0$ such that
\begin{align}
\|\cM f\|_p\leq C_{d,p}\|f\|_p.   \label{mf}
\end{align}
(ii) For every $f\in W^{1,p}_{loc}$, there is a constant $C_d>0$ such that for a.e. $x,y\in\mR^d$,
\begin{align}
|f(x)-f(y)|\leq C_{d}|x-y|\Big(\cM|\nabla f|(x)+\cM|\nabla f|(y)\Big).   \label{w11}
\end{align}
\el

The following multiplier theorem which gives $L^p$-boundness of operators can be found in \cite[Chapter IV, Section 3, Theorem 3]{St}.

\bl\label{mul}
Let $m(x)$ be a bounded function on $\mR^d$. Define the linear operator $T_m$ with domain $L^2(\mR^d)\cap L^p(\mR^d)$ by
the following relationship between Fourier transforms:
$$
\widehat{T_mf}(x)=m(x)\hat f(x).
$$
Suppose that $m(x)$ is of class $C^k$ in the complement of the origin of $\mR^d$, for every $\ell=(\ell_1,\ell_2,\cdots,\ell_d)$,
\begin{align}
|(\p_x)^{\ell}m(x)|\leq N|x|^{-|\ell|},\quad\forall|\ell|\leq k,  \label{mm}
\end{align}
where $k>\frac{d}{2}$ is an integer, $N$ is a positive constant, $|\ell|=\ell_1+\ell_2+\cdots+\ell_d$. Then, for every $1<p<\infty$,
there is a constant $N_p>0$ such that
$$
\|T_mf\|_p\leq N_p\|f\|_p,\quad\forall f\in L^p(\mR^d).
$$
\el

We will need the following Banach space version of the Calder\'on-Zygmund theorem, see \cite[Theorem 2.4]{Kr}.

\bl\label{cz}
Let $H$ be a Banach space. For each $(t,s)\in\mD=\{(t,s)\in\mR\times\mR: t\neq s\}$, a bounded operator $\kappa(t,s)$ from $H$ to $H$ satisfies
\begin{align*}
\|\nabla_s\kappa(t,s)\|_{L(H)}\leq N|t-s|^{-2},
\end{align*}
where $N>0$ is a constant, $\|\cdot\|_{L(H)}$ denotes the operator norm on $H$. Given $p\in(1,\infty)$, let $A: L^p(\mR,H)\rightarrow L^p(\mR,H)$
be a linear bounded operator with the form that: if $f$ is a bounded $H$-valued function and has compact support $\Gamma$, then for every $t\notin \Gamma$,
$$
Af(t)=\int_{-\infty}^{+\infty}\kappa(t,s)f(s)\dif s.
$$
Then, the operator $A$ is uniquely extendible to a bounded operator from $L^q(\mR,H)$ to $L^q(\mR,H)$ for any $q\in(1,p]$.
\el

Finally, we state the Khasminskii's estimate (see \cite[Lemma 13]{F-F-P-V}), which gives conditions ensuring the exponential integrability of a Markov process.

\bl\label{kha}
Let $(\Omega,\sF,(\mP_x)_{x\in\mR^{2d}}; (X_t)_{t\geq 0})$ be a family of $\mR^d$-valued time-homogenous
Markov process and $f$ be a non-negative measurable function on $\mR^d$. For any given $T>0$, if
\begin{align*}
\sup_{x\in\mR^{2d}}\mE_{x}\Bigg(\int_0^T\!\!f(X_t)\dif t\Bigg)=:c<1,
\end{align*}
where $\mE_x$ denotes the expectation with respect to $\mP_x$, then
\begin{align*}
\sup_{x\in\mR^{2d}}\mE_{x}\exp\Bigg\{\int_0^T\!f(X_t)\dif t\Bigg\}\leq \frac{1}{1-c}.
\end{align*}
\el

\section{Mixed type Kolmogorov equation in Bessel spaces}

Fix $T>0$ below and let $z=(x,y)\in\mR^{d_1+d_2}$. The functions $F$ and $G$ are given as in (\ref{sde1}). This section is devoted to the study of
the following partial integral-differential equation (PIDE) of Kolmogorov type over $[0,T]\times \mR^{d_1+d_2}$:
\begin{align}
\p_t u(t,z)=&\tfrac{1}{2}\Delta_yu(t,z)+\Delta^{\frac{\alpha}{2}}_xu(t,z)\no\\
&+F(t,z)\cdot\nabla_x u(t,z)+G(t,z)\cdot\nabla_y u(t,z)+f(t,z),\quad u(0,z)=0,   \label{pde1}
\end{align}
where $\alpha\in(1,2)$, $f:[0,T]\times\mR^{d_1+d_2}\rightarrow\mR^{d}$ is a measurable function.
To shorten the notation, for $p,q\in(1,\infty)$ and $\beta\in(0,2]$, we write
$$
\mL^p_{q}(T):=L^q\Big([0,T];L^p\big(\mR^{d_2}_y;L^p(\mR^{d_1}_x)\big)\Big)
$$
and
$$
\mH^{\beta,p}_{q}(T):=L^q\Big([0,T];L^p\big(\mR^{d_2}_y;H^{\beta}_p(\mR^{d_1}_x)\big)\Big),
$$
%and
%$$
%\mW^{\beta,p}_{q}(T):=L^q\Big([0,T];L^p\big(\mR^{d}_y;W^{\beta,p}(\mR^d_x)\big)\Big),
%$$
where we use $\mR^{d_1}_x, \mR^{d_2}_y$ to distinguish the space variables $x, y$.

\subsection{Equation (\ref{pde1}) without drift.} 
Let us first consider the following simpler PIDE on $[0,T]\times\mR^{d_1+d_2}$:
\begin{equation} \label{pde3}
\left\{ \begin{aligned}
&\p_t u(t,z)=\tfrac{1}{2}\Delta_yu(t,z)+\Delta^{\frac{\alpha}{2}}_xu(t,z)+f(t,z),\\
&u(0,z)=0.
\end{aligned} \right.
\end{equation}
For simplify, we write
$$
\sL_0\varphi(x,y):=\tfrac{1}{2}\Delta_y\varphi(x,y)+\Delta^{\frac{\alpha}{2}}_x\varphi(x,y),\quad \forall \varphi\in C_b^2(\mR^{d_1+d_2}).
$$
Consider the $\mR^{d_1+d_2}$-valued process
$$
U_t:=(L_t,W_t).
$$
where $L_t$ is the $d_1$-dimensional rotationally symmetric $\alpha$-stable process and $W_t$ is the $d_2$-dimensional Brownian motion.
Then, one can check easily that the infinitesimal operator of $U_t$ is just given by $\sL_0$.
Moreover, since we can always assume that $L_t$ and $W_t$ are independent, we have that, for any $z=(x,y)\in\mR^{d_1+d_2}$,
$$
\mP(U_t\leq z)=\mP(L_t\leq x,W_t\leq y)=\mP(L_t\leq x)\mP(W_t\leq y).
$$
Thus, the heat kernel $p(t,z)=p(t,x,y)$ of $U_t$ (which is also the fundamental solution of $\sL_0$) is given by
\begin{align}
p(t,x,y)=p^{(\alpha)}(t,x)p^{(2)}(t,y),  \label{pp}
\end{align}
where $p^{(\alpha)}(t,x)$ and $p^{(2)}(t,y)$ are the heat kernels of $L_t$ and $W_t$, respectively.

It is well known that
\begin{align}
p^{(2)}(t,y)=(\sqrt{2\pi}t)^{-d_2/2}\exp\Big\{-\frac{|y|^2}{2t}\Big\}  \label{h1}
\end{align}
and there exist constants $c_i>0$ ($i=1,2$) such that
\begin{align}
|\nabla^i_yp^{(2)}(t,y)|\leq c_it^{-\frac{d_2+i}{2}}\exp\Big\{-c_i\frac{|y|^2}{t}\Big\}.  \label{h2}
\end{align}
For $p^{(\alpha)}(t,x)$, we have the following estimates: there exist constents $c_0>1$ and $c_i>0$ ($i=1,2$) such that
\begin{align}
c_0^{-1}t\big(|x|+t^{\frac{1}{\alpha}}\big)^{-d_1-\alpha}\leq p^{(\alpha)}(t,x)\leq c_0t\big(|x|+t^{\frac{1}{\alpha}}\big)^{-d_1-\alpha}  \label{h3}
\end{align}
and
\begin{align}
|\nabla^i_xp^{(\alpha)}(t,x)|\leq c_it^{1-\frac{i}{\alpha}}\big(|x|+t^{\frac{1}{\alpha}}\big)^{-d_1-\alpha}.  \label{h4}
\end{align}
Denote by $\cT_t^{(2)}$ and $\cT^{(\alpha)}_t$ semigroup corresponding to $W_t$ and $L_t$, respectively.
Using (\ref{h1}), (\ref{h3}) and an easy computations yield that for every $p>1$, there exists a constant $c_{d_2}>0$ such that
\begin{align}
\|p^{(2)}(t,\cdot)\|_{L^{p}(\mR^{d_2}_y)}\leq c_{d_2} t^{-\frac{d_2(p-1)}{2p}}, \quad
\|p^{(\alpha)}(t,\cdot)\|_{L^{p}(\mR^{d_1}_x)}\leq c_{d_2}t^{-\frac{d_1(p-1)}{\alpha p}}.   \label{ppp}
\end{align}
Using H\"older's inequality, we have
\begin{align}
\|T_t^{(2)}f\|_{\infty}\leq \left(\int_{\mR^{d_2}}\big[ p^{(2)}(t,y)\big]^{p^*}\dif y\right)^{p^*}\cdot\|f\|_p\leq C_{d_2}t^{-\frac{d_2}{2p}}\|f\|_p,\label{b2}
\end{align}
where $C_{d_2}>0$, $p^*=\tfrac{p}{p-1}$ is the conjugate index of $p$. Similarly, we have
\begin{align}
\|T_t^{(\alpha)}f\|_{\infty}\leq C_{d_1}t^{-\frac{d_1}{\alpha p}}\|f\|_p.  \label{b1}
\end{align}
By (\ref{h2}) and (\ref{h4}), we have
\begin{align}
\|\nabla T_t^{(2)}f\|_{\infty}\leq C_{d_2}t^{-\frac{d_2+p}{2p}}\|f\|_p,\qquad\|\nabla T_t^{(\alpha)}f\|_{\infty}\leq C_{d_1}t^{-\frac{d_1+p}{\alpha p}}\|f\|_p.  \label{h5}
\end{align}

Let $\cT_t$ be the semigroup corresponding to $\sL_0$, i.e.,
$$
\cT_t\varphi(z)=\cT_t\varphi(x,y)=\int\!\!\!\!\int_{\mR^{d_1+d_2}}p(t,x-u,y-v)\varphi(u,v)\dif u\dif v,\quad\forall \varphi\in \sB_b(\mR^{d_1+d_2};\mR^d).
$$
We prove the following result.

\bl\label{lem}
For $p,q\in(1,\infty)$, $\beta\in[0,1]$ and $f\in \mH^{\beta,p}_{q}(T)$, there exists a unique strong solution
$u$ for PDE (\ref{pde3}) which is given by
\begin{align}
u(t,z)=\int_0^{t}\!\cT_{t-s}f(s,z) \dif s.  \label{uu}
\end{align}
Moreover, there is a constant $C=C(d,p,q,T)>0$ such that
\begin{align}
\|\p_tu\|_{\mH^{\beta,p}_{q}(T)}+\|\nabla^2_yu\|_{\mH^{\beta,p}_{q}(T)}+\|\nabla_yu\|_{\mH^{\frac{\alpha}{2}+\beta,p}_{q}(T)}+\|\Delta^{\frac{\alpha+\beta}{2}}_xu\|_{\mL^p_q(T)} \leq C\|f\|_{\mH^{\beta,p}_{q}(T)}.   \label{es1}
\end{align}
If we assume further that $p,q$ satisfying (\ref{index}), then we have for some $C_d>0$,
\begin{align}
\|\nabla_yu\|_{\infty}\leq C_dt^{\frac{1}{2}-\frac{1}{q}-\frac{d_2}{2p}-\frac{d_1}{\alpha p}}\|f\|_{\mL_q^p(T)},   \label{ess2}
\end{align}
and if $p,q$ satisfying
\begin{align}
\frac{d_2}{2p}+\frac{d_1}{\alpha p}+\frac{1}{q}<1-\frac{1}{\alpha},  \label{pqpq}
\end{align}
then we also have that
\begin{align}
\|\nabla_xu\|_{\infty}\leq C_dt^{1-\frac{1}{q}-\frac{d_2}{2p}-\frac{d_1}{\alpha p}-\frac{1}{\alpha}}\|f\|_{\mL_q^p(T)}.   \label{ess1}
\end{align}
\el

\begin{proof}
We first prove the desired result for $p=q$. Since $f(t,z)=f(t,x,y)$ can be approximated by the functions of the form
$f_1(t)f_2(x)f_3(y)$,
we only need to show that $u$ defined by (\ref{uu}) satisfies (\ref{pde3}) for $f_1\in L^p([0,T])$, $f_2\in H^{\beta}_p(\mR^{d_1}_x)$ and $f_3\in L^p(\mR^{d_2}_y)$.
In this case, we have
\begin{align}
u(t,z)&=\int_0^t\!\!\!\int_{\mR^{d_2}}\!\!\int_{\mR^{d_1}}p^{(\alpha)}(t-s,x-u)p^{(2)}(t-s,y-v)f_1(s)f_2(u)f_3(v)\dif u\dif v\dif s\no\\
&=\int_0^t\cT_{t-s}^{(\alpha)}f_2(x)\cT_{t-s}^{(2)} f_3(y)f_1(s)\dif s.\label{123}
\end{align}
It suffices to prove the estimate (\ref{es1}), the uniqueness follows by the standard method.

Note that by symmetric,
\begin{align*}
\Delta^{\frac{\beta}{2}}_xu(t,z)=\int_0^t\cT_{t-s}^{(\alpha)}\big(\Delta^{\frac{\beta}{2}}_xf_2\big)(x)\cT_{t-s}^{(2)} f_3(y)f_1(s)\dif s.
\end{align*}
Using Fourier transform, we can write
\begin{align*}
\cT_{t}^{(\alpha)}\big(\Delta^{\frac{\beta}{2}}_xf_2\big)(x)&=\Big(p^{(\alpha)}(t)*\big(\Delta^{\frac{\beta}{2}}f_2\big)\Big)(x)=\left(\widehat {p^{(\alpha)}(t)}\cdot \widehat{\big(\Delta^{\frac{\beta}{2}}_xf_2\big)}\right)^{\vee}(x)\\
&=\int_{\mR^{d_1}}\e^{ix\xi_2}\e^{-t|\xi_2|^\alpha}\widehat{\big(\Delta^{\frac{\beta}{2}}_xf_2\big)}(\xi_2)\dif \xi_2
\end{align*}
and
\begin{align*}
\cT_{t}^{(2)} f_3(y)&=\big(p^{(2)}(t)*f_3\big)(y)=\Big(\widehat {p^{(2)}(t)}\cdot \hat f_3\Big)^{\vee}(y)\\
&=\int_{\mR^{d_2}}\e^{iy\xi_3}\e^{-t|\xi_3|^2}\hat f_3(\xi_3)\dif \xi_3.
\end{align*}
Thus, letting $f_1(t)=0$ for $t<0$ and $t>T$, we have
\begin{align*}
u(t,z)&=\int_{-\infty}^{+\infty}\!\!\!\int_{\mR^{d_2}}\!\!\int_{\mR^{d_1}}\e^{i(x\xi_2+y\xi_3)}\e^{-(t-s)(|\xi_2|^{\alpha}+|\xi_3|^2)}
f_1(s)\widehat{\big(\Delta^{\frac{\beta}{2}}_xf_2\big)}(\xi_2)\hat f_3(\xi_3)\dif s\dif \xi_2\dif \xi_3\\
&=\int_{-\infty}^{+\infty}\!\!\!\int_{\mR^{d_2}}\!\!\int_{\mR^{d_1}}\e^{i(t\xi_1+x\xi_2+y\xi_3)}\frac{1}{i\xi_1+|\xi_2|^\alpha+|\xi_3|^2}
\bar f(\xi_1,\xi_2,\xi_3)\dif \xi_1\dif \xi_2\dif \xi_3,
\end{align*}
where
$$\hat f_\beta(\xi_1,\xi_2,\xi_3):=\hat f_1(\xi_1)\widehat{\big(\Delta^{\frac{\beta}{2}}_xf_2\big)}(\xi_2)\hat f_3(\xi_3).$$
As a result, we can get
\begin{align*}
\nabla^2_y\Delta^{\frac{\beta}{2}}_xu(t,z)=\int_{-\infty}^{+\infty}\!\!\!\int_{\mR^{d_2}}\!\!\int_{\mR^{d_1}}\e^{i(t\xi_1+x\xi_2+y\xi_3)}\frac{|\xi_3|^2}{i\xi_1+|\xi_2|^\alpha+|\xi_3|^2}
\hat f_\beta(\xi_1,\xi_2,\xi_3)\dif \xi_1\dif \xi_2\dif \xi_3,
\end{align*}
\begin{align*}
\Delta^{\frac{\alpha+\beta}{2}}_xu(t,z)=\int_{-\infty}^{+\infty}\!\!\!\int_{\mR^{d_2}}\!\!\int_{\mR^{d_1}}\e^{i(t\xi_1+x\xi_2+y\xi_3)} \frac{|\xi_2|^\alpha}{i\xi_1+|\xi_2|^\alpha+|\xi_3|^2}\hat f_\beta(\xi_1,\xi_2,\xi_3)\dif \xi_1\dif \xi_2\dif \xi_3,
\end{align*}
\begin{align*}
\p_t\Delta^{\frac{\beta}{2}}_xu(t,z)=\int_{-\infty}^{+\infty}\!\!\!\int_{\mR^{d_2}}\!\!\int_{\mR^{d_1}}\e^{i(t\xi_1+x\xi_2+y\xi_3)}\frac{i\xi_1}{i\xi_1+|\xi_2|^\alpha+|\xi_3|^2}
\hat f_\beta(\xi_1,\xi_2,\xi_3)\dif \xi_1\dif \xi_2\dif \xi_3
\end{align*}
and
\begin{align*}
\Delta^{\frac{\alpha/2+\beta}{2}}_x\nabla_yu(t,z)=\int_{-\infty}^{+\infty}\!\!\!\int_{\mR^{d_2}}\!\!\int_{\mR^{d_1}}\e^{i(t\xi_1+x\xi_2+y\xi_3)}
\frac{|\xi_2|^{\frac{\alpha}{2}}\cdot\xi_3}{i\xi_1+|\xi_2|^\alpha+|\xi_3|^2}\hat f_\beta(\xi_1,\xi_2,\xi_3)\dif \xi_1\dif \xi_2\dif \xi_3.
\end{align*}

Now, using Lemma \ref{mul} with
$$
m_1(\xi_1)=\frac{i\xi_1}{i\xi_1+|\xi_2|^\alpha+|\xi_3|^2},\,\,\,\, m_2(\xi_2)=\frac{|\xi_2|^\alpha}{i\xi_1+|\xi_2|^\alpha+|\xi_3|^2},
$$
$$
m_3(\xi_3)=\frac{|\xi_3|^2}{i\xi_1+|\xi_2|^\alpha+|\xi_3|^2},\,\,\,\, m_4(\xi_2,\xi_3)=\frac{|\xi_2|^{\alpha/2}\cdot\xi_3}{i\xi_1+|\xi_2|^\alpha+|\xi_3|^2},
$$
and checking as in \cite[Chapter IV, Section 3]{La-So-Ur} that (\ref{mm}) hold, for example,
\begin{align*}
|\xi_2|\cdot|\p_{\xi_2}m_4(\xi_2,\xi_3)|+|\xi_3|\cdot|\p_{\xi_3}m_4(\xi_2,\xi_3)|\leq N\frac{|\xi_2|^{\alpha/2}\cdot|\xi_3|}{c_0+|\xi_2|^\alpha+|\xi_3|^2}\leq N,
\end{align*}
we can conclude that
\begin{align*}
\|\p_tu\|_{\mH^{\beta,p}_{p}(T)}+\|\nabla^2_yu\|_{\mH^{\beta,p}_{p}(T)}&+\|\nabla_yu\|_{\mH^{\frac{\alpha}{2}+\beta,p}_{q}(T)}+\|\Delta^{\frac{\alpha}{2}}_xu\|_{\mH^{\beta,p}_{p}(T)}\\
&\leq C\|f_1(\Delta^{\frac{\beta}{2}}_xf_2)f_3\|_{\mL^p_p(T)}\leq C\|f\|_{\mH^{\beta,p}_{p}(T)}.
\end{align*}
This finished the proof when $p=q$.

For $p\neq q$, let us introduce the operators
$$
\kappa_1(t):=\Delta_y T_t^{(2)}T_t^{(\alpha)}
$$
and
$$
\kappa_2(t):=\Delta^{\frac{\alpha}{2}}_x T_t^{(2)}T_t^{(\alpha)}.
$$
Then, we have
$$
\Delta_yu(t,z)=\int_{-\infty}^{+\infty}\kappa_1(t-s)f(s)(x,y)\dif s
$$
and
$$
\Delta^{\frac{\alpha}{2}}_xu(t,z)=\int_{-\infty}^{+\infty}\kappa_2(t-s)f(s)(x,y)\dif s.
$$
It can be checked by (\ref{h2}) and (\ref{h4}) that
\begin{align*}
\|\p_s\kappa_1(t-s)f\|_{L^p(\mR^{d_1+d_2})}+\|\p_s\kappa_2(t-s)f\|_{L^p(\mR^{d_1+d_2})}\leq N(t-s)^{-2}\|f\|_{L^p(\mR^{d_1+d_2})},
\end{align*}
where $N$ is a positive constant.
Thus, Lemma \ref{cz} implies that $\Delta_y u$ and $\Delta^{\frac{\alpha}{2}}_xu$ are bounded operators from $L^q(\mR;L^p(\mR^{d_1+d_2}))$
into itself for $1<q\leq p$. By duality, we can conclude that $\Delta_y u$ and $\Delta^{\frac{\alpha}{2}}_xu$ are bounded in $L^q(\mR;L^p(\mR^{d_1+d_2}))$ for any $p,q\in(1,\infty)$. The same holds for $\Delta^{\frac{\alpha}{4}}\nabla_yu$. As for $\p_tu$, it follows by the equation (\ref{pde3}).

\vspace{2mm}
We proceed to prove the estimate (\ref{ess1}). Consider first that $u$ defined by
(\ref{123}). If $p,q\in(1,\infty)$ satisfy (\ref{index}), (\ref{b1}), (\ref{h5}) and H\"older's inequality yield
\begin{align*}
\|\nabla_xu\|_\infty&\leq\int_0^t\|\nabla\cT_{t-s}^{(\alpha)}f_2\|_{\infty}\|\cT_{t-s}^{(2)} f_3\|_\infty f_1(s)\dif s\\
&\leq C_d\left(\int_0^t(t-s)^{-(\frac{d_2}{2p}+\frac{d_1}{\alpha p}+\frac{1}{\alpha})q^*}\dif s\right)^{\frac{1}{q^*}}\cdot\|f\|_{\mL^q_p(T)},
\end{align*}
which in turn implies the desired results. The general case follows by standard approximation argument. The estimate (\ref{ess2}) follows by using (\ref{b2}), (\ref{h5})  and the same argument.
\end{proof}

\subsection{Equation (\ref{pde1}) with drift.}

Now, let us consider the general case of PIDE (\ref{pde1}).
We prove the following result.

\bt\label{th2}
Let $p,q\in(1,\infty)$ satisfy (\ref{pqpq}) and $T>0$ be small enough. Assume that $F,G\in \mL^p_{q}(T)$. Then, for every $f\in \mL^p_{q}(T)$
there exists a unique strong solution $u$ to PIDE (\ref{pde1}) with
\begin{align*}
\|\p_tu\|_{\mL^p_q(T)}+\|\nabla^2_yu\|_{\mL^p_q(T)}+\|\Delta^{\frac{\alpha}{4}}_x\nabla_yu\|_{\mL^p_q(T)}+\|\Delta^{\frac{\alpha}{2}}_xu\|_{\mL^p_q(T)} \leq C\|f\|_{\mL^q_p(T)},
\end{align*}
where $C=C(d,\|F\|_{\mL^p_q(T)},\|G\|_{\mL^p_q(T)},T)$ is a positive constant. Moreover, estimates (\ref{ess2}) and (\ref{ess1}) hold.
\et

\begin{proof}
It is now standard to show that $u$ is a solution to PIDE (\ref{pde1}) if and only if $u$ satisfies the following integral equation (see \cite[Proposition 3.5]{Zh00}):
\begin{align}
u(t,z)=\int_0^t\cT_{t-s}\big(f+F\cdot\nabla_xu+G\cdot\nabla_yu\big)(s,z)\dif s.  \label{inte}
\end{align}
We seek a solution $u$ of (\ref{inte}) by the Picard iteration. Let $u_0\equiv0$ and for $n\in \mN$, define $u_n$ recursively by
\begin{align}
u_n(t,z):=\int_0^t\cT_{t-s}\big(f+F\cdot\nabla_xu_{n-1}+G\cdot\nabla_yu_{n-1}\big)(s,z)\dif s.  \label{nnn}
\end{align}
Then, by Lemma \ref{lem},
$$
u_1(t,z)=\int_0^t\cT_{t-s}f(s,z)\dif s$$
belong to $\mH^{\alpha,p}_q(T)\cap L^q\Big([0,T];W^{2,p}\big(\mR^{d_2}_y;L^p(\mR^{d_1}_x)\big)\Big)$.
Since $p, q$ satisfy (\ref{index}), estimates (\ref{ess1})-(\ref{ess2}) yield that $\nabla_x u_1$ and $\nabla_y u_1$ are all bounded. We then have
$$
F\cdot\nabla_xu_1,\,\, G\cdot\nabla_y u_1\in \mL^p_q(T).
$$
This means that $u_2$ is well defined, and so on.

By induction, we prove the following claim: there exist constants $C_1>0$ and $\delta\in(0,1)$ such that, for any $n\in\mN$,
\begin{align}
&\|\p_tu_n\|_{\mL^p_q(T)}+\|\nabla^2_yu_n\|_{\mL^p_q(T)}+\|\nabla_y\Delta^{\frac{\alpha}{4}}_xu_n\|_{\mL^p_q(T)}
+\|\Delta^{\frac{\alpha}{2}}_xu_n\|_{\mL^p_q(T)}\no\\
\leq &C_1\Big(1+\sum_{k=1}^{n-1}\delta^k\Big)\|f\|_{\mL^p_q(T)}  \label{nn}
\end{align}
and
\begin{align}
\|\nabla_xu_n\|_{\infty}+\|\nabla_y u_n\|_{\infty}\leq 2C_1T^{1-\frac{1}{\alpha}-\frac{1}{q}-\frac{d_2}{2p}-\frac{d_1}{\alpha p}}\Big(1+\sum_{k=1}^{n-1}\delta^k\Big)\|f\|_{\mL^p_q(T)}.   \label{nn2}
\end{align}
{\it Proof of Claim}:
Choosing $T<1$, (\ref{es1})-(\ref{ess2}) imply that (\ref{nn}) and (\ref{nn2}) holds for $n=1$. Assume that (\ref{nn})-(\ref{nn2})
hold for $n\geq 1$. We proceed to show that they are still true for $n+1$. In fact, we have
\begin{align*}
&\|\p_tu_{n+1}\|_{\mL^p_q(T)}+\|\nabla^2_yu_{n+1}\|_{\mL^p_q(T)}+\|\Delta^{\frac{\alpha}{4}}_x\nabla_yu_n\|_{\mL^p_q(T)}+\|\Delta^{\frac{\alpha}{2}}_xu_{n+1}\|_{\mL^p_q(T)}\\
\leq& C_1\|f-F\cdot\nabla_xu_n-G\cdot\nabla_yu_n\|_{\mL^p_q(T)}\\
\leq& C_1\Big(\|f\|_{\mL^p_q(T)}+\|F\cdot\nabla_x u_n\|_{\mL^p_q(T)}+\|G\cdot\nabla_yu_n\|_{\mL^p_q(T)}\Big)\\
\leq& C_1\Big(\|f\|_{\mL^p_q(T)}+\|F\|_{\mL^p_q(T)}\|\nabla_x u_n\|_\infty+\|G\|_{\mL^p_q(T)}\|\nabla_yu_n\|_\infty\Big).
\end{align*}
By the induction assumption (\ref{nn2}), we can get
\begin{align*}
&\|F\|_{\mL^p_q(T)}\|\nabla_x u_n\|_\infty+\|G\|_{\mL^p_q(T)}\|\nabla_yu_n\|_\infty\\
&\leq 2C_1T^{1-\frac{1}{\alpha}-\frac{1}{q}-\frac{d_2}{2p}-\frac{d_1}{\alpha p}}\Big(\|F\|_{\mL^p_q(T)}+\|G\|_{\mL^p_q(T)}\Big)\cdot\Big(1+\sum_{k=1}^{n-1}\delta^k\Big)\|f\|_{\mL^p_q(T)}.
\end{align*}
Choosing $T$ small enough such that
$$
2C_1T^{1-\frac{1}{\alpha}-\frac{1}{q}-\frac{d_2}{2p}-\frac{d_1}{\alpha p}}\Big(\|F\|_{\mL^p_q(T)}+\|G\|_{\mL^p_q(T)}\Big)<\delta,
$$
we get the desired result (\ref{nn}). As for (\ref{nn2}), it follows by the same argument.

Thus, we obtain that there is a constant $C_2>0$ such that
\begin{align*}
\sup_{n\in\mN}\Big(\|\p_tu_n\|_{\mL^p_q(T)}+\|\nabla^2_yu_n\|_{\mL^p_q(T)}+\|\Delta^{\frac{\alpha}{4}}_x\nabla_yu_n\|_{\mL^p_q(T)}+\|\Delta^{\frac{\alpha}{2}}_xu_n
\|_{\mL^p_q(T)}\Big)\leq C_2\|f\|_{\mL^p_q(T)}
\end{align*}
and
\begin{align*}
\sup_{n\in\mN}\Big(\|u_n\|_{\infty}+\|\nabla u_n\|_{\infty}\Big)\leq C_2\|f\|_{\mL^p_q(T)}.
\end{align*}
We may also argue as before to obtain that for $n>m$,
\begin{align*}
\|\p_tu_n\|_{\mL^p_q(T)}+\|\nabla^2_yu_n\|_{\mL^p_q(T)}+\|\Delta^{\frac{\alpha}{2}}_xu_n\|_{\mL^p_q(T)}\leq C_1\left(\sum_{k=m+1}^n\delta^k\right)\|f\|_{\mL^q_p(S,T)}.
\end{align*}
Hence, $\{u_n\}_{n\geq 1}$ is a Cauchy sequence in
$$
W^{1,q}\Big([0,T];L^p\big(\mR^{d_2}_y;L^p(\mR^{d_1}_x)\big)\Big)\cap L^q\Big([0,T];W^{2,p}\big(\mR^{d_2}_y;L^p(\mR^{d_1}_x)\big)\Big)\cap \mH^{\alpha,p}_q(T).
$$
As a result, there exists a function $u$ such that
$$
\lim_{n\rightarrow\infty}\Big(\|\nabla^2 (u_n-u)\|_{\mL^q_p(0,T)}+\|\p_t (u_n-u)\|_{\mL^q_p(0,T)}\Big)=0.
$$
Taking limits for both sides of (\ref{nnn}), we obtain the existence of a solution. The proof is finished.
\end{proof}

The higher regularity in $x$ is not obvious as usual since we do not assume that $F$ is bounded. As a result, we need to make full use of Sobolev embedding theorems.
\bt\label{th3}
Let $p,q\in(1,\infty)$ satisfy (\ref{index}) and $\beta\in(0,1]$ satisfy (\ref{index2}). Assume that $F,G\in \mH^{\beta,p}_{q}(T)$.
Then, for every $f\in \mH^{\beta,p}_{q}(T)$, there exists a unique strong solution
$u$ to PIDE (\ref{pde1}) with estimates (\ref{es1})-(\ref{ess2}) hold.
\et

\begin{proof}
We only show the priori estimate, then one can follow the same procedure as in the proof of Theorem \ref{th2}. For this purpose,
we assume that there exists a function $u$ satisfying the PIDE (\ref{pde1}) and (\ref{es1})-(\ref{ess2}). Then, by Lemma \ref{lem} we can derive that
\begin{align*}
&\|\p_tu\|_{\mH^{\beta,p}_q(T)}+\|\Delta^{\frac{\alpha}{2}}_xu\|_{\mH^{\beta,p}_q(T)}+\|\Delta^{\frac{\alpha}{4}}_x\nabla_yu\|_{\mL^p_q(T)}+\|\nabla^2_yu\|_{\mH^{\beta,p}_q(T)}\\
&\leq C\Big(\|f\|_{\mH^{\beta,p}_{q}(T)}+\|F\cdot\nabla_x u\|_{\mH^{\beta,p}_q(T)}+\|G\cdot\nabla_y u\|_{\mH^{\beta,p}_q(T)}\Big).
\end{align*}
We proceed to show that terms on the right hand of the above inequality are well defined. Obviously, we have by assumption that
$$
\|f\|_{\mH^{\beta,p}_{q}(T)}<\infty.
$$
Thanks to the assumption that $p,q$ satisfy (\ref{index}) and $\beta$ satisfy (\ref{index2}), and by Sobolev embedding (\ref{emb2}), we have
$f\in \mL^{\hat q}_{\hat p}(T)$ with $\hat p,\hat q$ satisfying (\ref{pqpq}). As a result, (\ref{ess1}) is true and
by the definition, we can write
\begin{align}
\|F\cdot\nabla_x u\|_{\mH^{\beta,p}_q(T)}\leq \|F\|_{\mL^{p}_q(T)}\|\nabla_xu\|_\infty+\|F\|_{\mH^{\beta,p}_q(T)}
\|\nabla_xu\|_\infty+\|F\cdot\Delta^{\frac{\beta}{2}}_x\nabla_xu\|_{\mL^{p}_q(T)}.\label{77}
\end{align}
The first two terms can be controlled easily due to the assumption and (\ref{ess1}). As for the third term,
we'll discuss in three different cases.

\vspace{3mm}
{\it Case 1: $\beta<\frac{d_1}{p}$}. On one hand,
we have by the Sobolev embedding (\ref{emb2}) that for every $0<\eps<\beta$,
$$
F\in \mH^{\beta,p}_q(T)\subseteq\mH^{\beta-\eps,p}_q(T)\hookrightarrow L^q\Big([0,T];L^p\big(\mR^{d_2}_y;L^{p^*}(\mR^{d_1}_x)\big)\Big)
$$
with $p^*$ be the Sobolev conjugate index of $p$ given by
$$
p^*=\frac{d_1}{d_1/p-\beta+\eps}.
$$
On the other hand, by (\ref{index}) and (\ref{index2}) we can choose $\eps$ small enough such that
$$
\alpha-1>\frac{d_1}{p}-\beta+\eps.
$$
(\ref{re}), (\ref{emb2}) and (\ref{emb}) imply that
\begin{align*}
&\Delta^{\frac{\beta}{2}}_x\nabla_x u\in \mH^{\alpha-1,p}_q(T)\cap W^{1,q}\Big([0,T];L^p\big(\mR^{d_2}_y;L^{p}(\mR^{d_1}_x)\big)\Big)\cap L^q\Big([0,T];W^{2,p}\big(\mR^{d_2}_y;L^{p}(\mR^{d_1}_x)\big)\Big)\\
&\hookrightarrow\mH^{\frac{d}{p}-\beta+\eps,p}_q(T)\cap W^{1,q}\Big([0,T];L^p\big(\mR^{d_2}_y;L^{p}(\mR^{d_1}_x)\big)\Big)\cap L^q\Big([0,T];W^{2,p}\big(\mR^{d_2}_y;L^{p}(\mR^{d_1}_x)\big)\Big)\\
&\hookrightarrow L^\infty\Big([0,T];L^\infty\big(\mR^{d_2}_y;L^{r}(\mR^{d_1}_x)\big)\Big)
\end{align*}
with
$$
r=\frac{d_1}{d_1/p-d_1/p+\beta-\eps}=\frac{d_1}{\beta-\eps}.
$$
It is clear that
$$
\frac{1}{p}=\frac{1}{p^*}+\frac{1}{r}.
$$
By H\"older's inequality, it holds
\begin{align*}
\|F\cdot\Delta^{\frac{\beta}{2}}_x\nabla_x u\|_{\mL^{p}_q(T)}&\leq \|F\|_{L^q([0,T];L^p(\mR^{d_2}_y;L^{p^*}(\mR^{d_1}_x)))}\|\Delta^{\frac{\beta}{2}}_x\nabla_x u\|_{L^\infty([0,T];L^\infty(\mR^{d_2}_y;L^{r}(\mR^{d_1}_x)))}\\
&\leq \|F\|_{\mH^{\beta,p}_q(T)}\|u\|_{\mH^{\alpha+\beta,p}_q(T)}<\infty.
\end{align*}
Consequently, every terms in the right of (\ref{77}) are finite.

\vspace{3mm}
{\it Case 2: $\beta=\frac{d_1}{p}$}. In this case, we can still use (\ref{re}) to get that for $\eps>0$ small enough,
$$
F\in \mH^{\beta,p}_q(T)\subseteq\mH^{\beta-\eps,p}_q(T).
$$
Thus, we may argue entirely the same as the above to get the desired results.

\vspace{3mm}
{\it Case 3: $\beta>\frac{d_1}{p}$}. By the Sobolev embedding (\ref{emb}), $F$ is in fact H\"older continuity in $x$ and we have
$$
F\in \mH^{\beta,p}_q(T)\hookrightarrow L^q\Big([0,T];L^p\big(\mR^{d_2}_y;L^{\infty}(\mR^{d_1}_x)\big)\Big).
$$
Meanwhile, we also have
\begin{align*}
\Delta^{\frac{\beta}{2}}_x\nabla_x u&\in \mH^{\alpha-1,p}_q(T)\cap W^{1,q}\Big([0,T];L^p\big(\mR^{d_2}_y;L^{p}(\mR^{d_1}_x)\big)\Big)\cap L^q\Big([0,T];W^{2,p}\big(\mR^{d_2}_y;L^{p}(\mR^{d_1}_x)\big)\Big)\\
&\hookrightarrow L^\infty\Big([0,T];L^\infty\big(\mR^{d_2}_y;L^{p}(\mR^{d_1}_x)\big)\Big).
\end{align*}
Hence, it holds that
\begin{align*}
\|F\cdot\Delta^{\frac{\beta}{2}}_x\nabla_x u\|_{\mL^p_q(T)}&\leq\|F\|_{L^q([0,T];L^p(\mR^{d_2}_y;L^{\infty}(\mR^{d_1}_x)))}\|\Delta^{\frac{\beta}{2}}_x\nabla_x u
\|_{L^\infty([0,T];L^\infty(\mR^{d_2}_y;L^{p}(\mR^{d_1}_x)))}\\
&\leq \|F\|_{\mH^{\beta,p}_q(T)}\|u\|_{\mH^{\alpha+\beta,p}_q(T)}<\infty.
\end{align*}
Combing the above computations, we get that
$$
\|F\cdot\nabla_x u\|_{\mH^{\beta,p}_q(T)}\leq \|F\|_{\mH^{\beta,p}_q(T)}\|u\|_{\mH^{\alpha+\beta,p}_q(T)}<\infty.
$$
Following the same argument, we can show that
$$
\|G\cdot\nabla_y u_n\|_{\mH^{\beta,p}_q(T)}\leq \|G\|_{\mH^{\beta,p}_q(T)}\|u\|_{\mH^{\alpha+\beta,p}_q(T)}<\infty.
$$
The proof can be finished.
\end{proof}

\section{Pathwise uniqueness of strong solutions}

This section will be divided into three parts: firstly, we prove the Krylov estimate for the strong solutions of SDE (\ref{sde1}).
Then, we perform the Zvonkin transformation by using the results in Section 3. Finally, we give the proof of our main results.

\subsection{Krylov estimate}

The method we use here is the Girsanov transformation. First, we prove the Krylov estimate $U_t$.
\bl\label{kry1}
Let $f\in\mL^q_p(T)$ be a non-negative function with $p, q$ satisfying
\begin{align}
\frac{d_2}{2p}+\frac{d_1}{\alpha p}+\frac{1}{q}<1.    \label{pq}
\end{align}
Then, there exists a constant $C=C(d_1,d_2,T,p,q)$ such that
$$
\mE\left(\int_0^tf(s,U_s)\dif s\right)\leq C\|f\|_{\mL^p_q(T)}.
$$
\el
\begin{proof}
In fact, by (\ref{pp}) and H\"older's inequality we have
\begin{align*}
\mE\left(\int_0^tf(s,U_s)\dif s\right)&=\int_0^t\!\!\!\int\!\!\!\!\int_{\mR^{d_1+d_2}}f(s,x,y)p^{(\alpha)}(s,x)p^{(2)}(s,y)\dif x\dif y\dif s\\
&\leq\int_0^t\!\!\!\int_{\mR^{d_1}}\|f(s,x,\cdot)\|_{L^p(\mR^{d_2}_y)}p^{(\alpha)}(s,x)\dif x\cdot\|p^{(2)}(s,\cdot)\|_{L^{p^*}(\mR^{d_2}_y)}\dif s\\
&\leq \int_0^t\|f(s,\cdot,\cdot)\|_{L^p(\mR^{d_1}_x;L^p(\mR^{d_2}_y))}\|p^{(\alpha)}(s,\cdot)\|_{L^{p^*}(\mR^{d_1}_x)}\cdot\|p^{(2)}(s,\cdot)\|_{L^{p^*}(\mR^{d_2}_y)}\dif s\\
&\leq \|f\|_{\mL^p_q(t)}\left(\int_0^t\|p^{(\alpha)}(s,\cdot)\|_{L^{p^*}(\mR^{d_1}_x)}^{q^*}\cdot\|p^{(2)}(s,\cdot)\|_{L^{p^*}(\mR^{d_2}_y)}^{q^*}\dif s\right)^{\frac{1}{q^*}},
\end{align*}
where $p^*$, $q^*$ are the conjugate index of $p$ and $q$, respectively.
Hence, the assumption (\ref{pq}) and (\ref{ppp}) imply that
$$
\int_0^ts^{-(\frac{d_1}{\alpha p}+\frac{d_2}{2 p})q^*}\dif s=\int_0^ts^{-(\frac{d_1}{\alpha p}+\frac{d_2}{2 p})\frac{q}{q-1}}\dif s<\infty.
$$
The proof is finished.
\end{proof}

As a direct consequence of Lemma \ref{kha}, we have the following result.

\bl
For every $T>0$ and non-negative function $f\in\mL^q_p(T)$ with $p, q$ satisfying (\ref{pq}), we have
\begin{align}
\sup_{z\in\mR^{d_1+d_2}}\mE^z\exp\bigg\{\int_0^Tf(s,U_s)\dif s\bigg\}\leq C<\infty,    \label{exes}
\end{align}
where $C>0$ is a constant depending on $p,q,d_1,d_2,T$ and $\|f\|_{\mL^q_p(T)}$.
\el

\begin{proof}
By the assumption, we can always chose a $r>1$ such that
$$
\frac{d_2r}{2p}+\frac{d_1r}{\alpha p}+\frac{r}{q}<1.
$$
Meanwhile, by Young's inequality, for any $\eps>0$, there exists a $C_\eps$ such that
$$
f\leq \eps|f|^r+C_\eps.
$$
In view of (\ref{kry1}), we can choose $\eps$ small enough such that
\begin{align}
\mE\left(\int_0^t\eps f^r(s,U_s)\dif s\right)\leq \eps C\|f^r\|_{\mL^q_p(T)}<1.  \label{fr}
\end{align}
Thus, we can deduce by Lemma \ref{kha} that
\begin{align*}
\mE\exp\bigg\{\int_0^Tf(s,U_s)\dif s\bigg\}\leq \e^{C_\eps T}\cdot\mE\exp\bigg\{\int_0^T\eps f^r(s,U_s)\dif s\bigg\}<\infty.
\end{align*}
The proof is finished.
\end{proof}

Now, using the Girsanov theorem, we show the Krylov estimate for
$(X_t, Y_t)$. For this purpose, Notice that
$$
U_t=QL_t+RW_t.
$$
We prove the following important result.

\bl\label{weak}
Suppose that $F, G\in\mL^q_p(T)$ with $p, q$ satisfying (\ref{index}). Then, for each $z\in\mR^{d_1+d_2}$ there exists a weak solution $Z_t$
to SDE (\ref{sde2}). Moreover, for every non-negative function $f\in\mL^q_p(T)$ with $p, q$ satisfying (\ref{index}),
\begin{align}
\sup_{z\in\mR^{d_1+d_2}}\mE\left(\int_0^Tf(s,Z_s)\dif s\right)\leq C_1\|f\|_{\mL^p_q(T)},    \label{kry2}
\end{align}
where $C_1$ is a constant depending on $p,q,d_1,d_2,T$ and $\|F\|_{\mL^q_p(T)}$, $\|G\|_{\mL^q_p(T)}$. In particular, we have
\begin{align}
\sup_{z\in\mR^{d_1+d_2}}\mE\exp\bigg\{\int_0^Tf(s,Z_s)\dif s\bigg\}\leq C_2<\infty,    \label{exes2}
\end{align}
where $C_2$ is a constant depending on $p,q,d_1,d_2,T$, $\|F\|_{\mL^q_p(T)}$, $\|G\|_{\mL^q_p(T)}$ and $\|f\|_{\mL^q_p(T)}$.
\el

\begin{proof}
We argue similarly as in the proof of Theorem 15 in \cite{F-F-P-V}.
Let $(\Omega,\sF,\mP)$ be the probability space on which Brownian motion $W_t$ and L\'evy process $L_t$ are defined.
Define the process
$$
H_t:=RW_t-\int_0^tRF(s,U_s)\dif s-\int_0^tQG(s,U_s)\dif s
$$
and set
$$
\varphi_t:=\exp\Bigg\{\int_0^t\!\langle F(s,U_s)+G(s,U_s),\dif W_s\rangle-\frac{1}{2}\int_0^t\!|B(s,U_s)|^2\dif s\Bigg\}.
$$
The assumption (\ref{index}) yields that
$$
|B|^2\in \mL^{q'}_{p'}(T)\quad\text{with}\quad q'=q/2,\,\,p=p/2.
$$
One can check that
$$
\frac{d}{2p'}+\frac{d}{\alpha p'}+\frac{1}{q'}<2-\frac{2}{\alpha}<1.
$$
Thus, using (\ref{exes}) for $f=|B|^2$, the Novikov condition ensures that the process $\varphi_T$ is a martingale.
Set the new probability measure by
$$
\frac{\dif \mQ}{\dif \mP}=\varphi_T.
$$
Then, by the Girsanov theorem, we know that $H_t$ is a Brownian motion under the new measure $\mQ$ with matrix given by
\begin{equation*}
\left(
  \begin{array}{l}
    0\,\, 0\\
    0\,\, \mI_d \\
  \end{array}
\right).
\end{equation*}
That is, $H_t$ has the same distribution under $\mQ$ as $RW_t$ under $\mP$.
As a result, we have
\begin{align*}
\dif U_t&=R\dif W_t+Q\dif L_t=\dif H_t+Q\dif L_t+RF(t,U_t)\dif t+QG(t,U_t)\dif t\\
&=B(t,U_t)\dif t+Q\dif L_t+\dif H_t,
\end{align*}
which means that $U_t$ is a solution under the probability measure $\mQ$.

Next, we proceed to prove the estimate (\ref{kry2}). We have by H\"older's inequality that, for some $r>1$ and $1/r+1/r'=1$,
\begin{align*}
\mE^{\mP}\left(\int_0^Tf(s,Z_s)\dif s\right)=\mE^{\mQ}\left(\int_0^Tf(s,U_s)\dif s\right)\leq \bigg[\mE^{\mP}\int_0^Tf^r(s,U_s)\dif s\bigg]^{1/r}
\Big[\mE^{\mP}\big(\varphi_T^{r'}\big)\Big]^{1/r'}
\end{align*}
The same argument as in (\ref{fr}), we can choose $r>1$ small enough such that
$$
\mE^{\mP}\left(\int_0^Tf^r(s,U_s)\dif s\right)<\infty.
$$
Then, we can write
\begin{align*}
\varphi_t^{r'}=\exp\Bigg\{\int_0^t r'\langle F(s,U_s)&+G(s,U_s),\dif W_s\rangle-\frac{1}{2}\int_0^t|r'B(s,U_s)|^2\dif s\\
&+\frac{r'(r'-1)}{2}\int_0^t|B(s,U_s)|^2\dif s\Bigg\}.
\end{align*}
Novikov condition and (\ref{exes}) ensures that $\varphi_T^{r'}$ has finite expectation, which in turn implies that (\ref{kry2}) is true. As for (\ref{exes2}),
it follows by the same way as in the proof of (\ref{exes}).
\end{proof}

\subsection{Zvonkin transformation}

Now, let us fix $T>0$ small enough, and consider the following PIDE in $[0,T]\times \mR^{d_1+d_2}$:
\begin{equation} \label{pde}
\left\{ \begin{aligned}
&\p_t u_1=\sL_0u_1+F\cdot\nabla_x u_1+G\cdot\nabla_y u_1+F,\\
&\p_t  u_2=\sL_0u_2+F\cdot\nabla_x u_2+G\cdot\nabla_y u_2+G.
\end{aligned} \right.
\end{equation}
According to Theorem \ref{th3}, there exist two functions $u_1, u_2$ which satisfy (\ref{pde}) and (\ref{es1}) holds. If we set
\begin{equation*}
U(t,z)=\left(
  \begin{array}{c}
    u_1(T-t,z) \\
    u_2(T-t,z) \\
  \end{array}
\right),
\end{equation*}
then the PIDE (\ref{pde}) can be written as
\begin{align}
\p_t U+\sL_0 U+\nabla U\cdot B+B=0, \label{pde2}
\end{align}
where
\begin{equation*}
\nabla U=\left(
  \begin{array}{c}
    \nabla u_1 \\
    \nabla u_2 \\
  \end{array}
\right)
=\left(
  \begin{array}{cc}
    \nabla_x u_1 & \nabla_y u_1\\
    \nabla_x u_2 &  \nabla_y u_2\\
  \end{array}
\right).
\end{equation*}

Let us define
\begin{align}
\Phi_t(z):=z+U(t,z). \label{phi}
\end{align}
Choosing $T$ small enough, (\ref{ess1}) and (\ref{ess2}) imply
$$
\frac{1}{2}|z_1-z_2|\leq\big|\Phi_t(z_1)-\Phi_t(z_2)\big|\leq \frac{3}{2}|z_1-z_2|,
$$
which implies that the map $z\to\Phi_t(z)$ forms a $C^1$-diffeomorphism and
\begin{align}
\frac{1}{2}\leq \|\nabla\Phi\|_{\infty},\|\nabla\Phi^{-1}\|_{\infty}\leq 2,   \label{upd}
\end{align}
where $\Phi^{-1}_t(\cdot)$ is the inverse function of $\Phi_t(\cdot)$. We prove the following Zvonkin's transformation.
%The main point is that since $b$ is un-bounded, we should be very careful with the index in order to use the Krylov estimate.

\bl\label{zvon}
Let $\Phi_t(z)$ be defined as the above and $(X_t,Y_t)$ solve SDE (\ref{sde1}). Then, $\hat Z_t:=\Phi_t(X_t,Y_t)$ satisfies the following SDE:
\begin{align*}
\hat Z_t&=\Phi_0(z)+\int_0^t\nabla_y \Phi_s\big(\Phi_s^{-1}(\hat Z_s)\big)\dif W_s+\int_0^t\tilde B\big(s,\Phi_s^{-1}(\hat Z_s)\big)\dif s\\
&\quad+\int_0^t\!\!\!\int_{|v|\leq 1}\Big[\Phi_s\big(\Phi_s^{-1}(\hat Z_s)+Qv\big)-\hat Z_s\Big]\tilde N(\dif v,\dif s)\\
&\quad+\int_0^t\!\!\!\int_{|v|> 1}\Big[\Phi_s\big(\Phi_s^{-1}(\hat Z_s)+Qv\big)-\hat Z_s\Big] N(\dif v,\dif s).
\end{align*}
where
$$
\tilde B(s,z):=-\int_{|v|>1}\Big[U\big(s,\Phi_s^{-1}(z)+Qv\big)-U\big(s,\Phi_s^{-1}(z)\big)\Big]\nu(\dif v).
$$
\el
\begin{proof}
Let $\rho$ be a non-negative smooth function on $\mR^{d_1+d_2+1}$ with support in $\{x\in\mR^{d_1+d_2+1}: |x|\leq 1\}$ and $\int_{\mR^{d_1+d_2+1}}\rho(t,y)
\dif t\dif y=1$. Set $\rho_n(t,y):=n^{d+1}\rho(nt,ny)$, and extend $U(t)$ to $\mR$ by setting $U(t,\cdot)=0$ for $t\geq T$ and $U(t,\cdot)=U(0,\cdot)$ for $t\leq 0$. Define
\begin{align*}
U_n(t,z):=\int_{\mR^{d_1+d_2+1}}U(s,y)\rho_n(t-s,z-y)\dif s\dif y,
\end{align*}
and set
\begin{align}
B_n(t,z):=-\p_t U_n(t,z)-\sL_0 U_n(t,z)-\nabla U_n(t,z)\cdot B(t,z). \label{eqq}
\end{align}
Thus, by (\ref{es1}), (\ref{pde2}) and the properties of convolution we get
\begin{align*}
\|B_n-B\|_{\mL^p_q(T)}&\leq \|\p_t (U_n-U)\|_{\mL^{p}_{q}(T)}+\|\nabla^2_y(U_n-U)\|_{\mL^{p}_{q}(T)}+\|\Delta^{\frac{\alpha}{2}}_x(U_n-U)\|_{\mL^{p}_{q}(T)}\\
&\quad+\|\nabla (U_n-U)\|_{\infty}\|B\|_{\mL^p_q(T)}\rightarrow0,\quad \text{as}\,\,n\rightarrow\infty.
\end{align*}
By Ito's formula, we have
\begin{align*}
U_n(t,Z_t)=&U_n(0,z)+\int_0^t\!\Big(\p_sU_n+\sL_0 U_n+\nabla U_n\cdot B\Big)(s,Z_s)\dif s+\int_0^t\nabla_y U_n(s,Z_s)\dif W_s\\
&-\int_0^t\!\!\!\int_{|v|>1}\Big[U_n\big(s,X_s+v,Y_s\big)-U_n(s,X_s,Y_s)\Big]\nu(\dif v)\dif s\\
&+\int_0^t\!\!\!\int_{|v|> 1}\Big[U_n\big(s,X_{s-}+v,Y_s\big)-U_n(s,X_{s-},Y_s)\Big] N(\dif v, \dif s)\\
&+\int_0^t\!\!\!\int_{|v|\leq 1}\Big[U_n\big(s,X_{s-}+v,Y_s\big)-U_n(s,X_{s-},Y_s)\Big]\tilde N(\dif v, \dif s).
\end{align*}
Set
$$
\Phi_t^n(z):=z+U_n(t,z).
$$
Adding the above equation with SDE \eqref{sde1}, by \eqref{eqq} and noticing that
$$
\Phi^n_t(z+Qv)-\Phi_t^n(z)=U_n(t,z+Qv)-U_n(t,z)+Qv,
$$
we obtain
\begin{align*}
\hat Z^n_t&:=\Phi_t^n(Z_t)=\Phi_0^n(z)+\int_0^t\!\Big(B(s,Z_s)-B_n(s,Z_s)\Big)\dif s+\int_0^t\nabla_y \Phi^n_s(\hat Z_s)\dif W_s\\
&\quad-\int_0^t\!\!\!\int_{|v|>1}\Big[U_n\big(s,Z_s+Qv\big)-U_n(s,Z_s)\Big]\nu(\dif v)\dif s\\
&\quad+\int_0^t\!\!\!\int_{|v|> 1}\Big[\Phi_t^n\big(Z_{s-}+Qv\big)-\Phi_t^n(Z_{s-})\Big] N(\dif v, \dif s)\\
&\quad+\int_0^t\!\!\!\int_{|v|\leq 1}\Big[\Phi_t^n\big(Z_{s-}+Qv\big)-\Phi_t^n(Z_{s-})\Big]\tilde N(\dif v, \dif s)\\
&=:\Phi_n(z)+\cQ_1+\cQ_2+\cQ_3+\cQ_4.
\end{align*}
Now we are going to take limit for the above equality. First of all, it is easy to see that
$$
\lim_{n\rightarrow\infty}\hat Z^n_t=\Phi_t(Z_t)=\hat Z_t.
$$
By the dominated convergence theorem, we also have
\begin{align*}
\cQ_2+\cQ_3\rightarrow&-\int_0^t\!\!\!\int_{|v|>1}\Big[U\big(s,Z_s+Qv\big)-U(s,Z_s)\Big]\nu(\dif v)\dif s\\
&+\int_0^t\!\!\!\int_{|v|> 1}\Big[\Phi_t\big(Z_{s-}+Qv\big)-\Phi_t(Z_{s-})\Big] N(\dif v, \dif s).
\end{align*}
As for $\cQ_4$, (\ref{upd}) and the dominated convergence theorem imply that
\begin{align*}
&\mE\bigg|\int_0^t\!\!\!\int_{|v|\leq1}\!\!\Big[\Phi_s^n(Z_{s-}\!+\!Qv)\!-\!\Phi_s^n(Z_{s-})\!-\!\Phi(Z_{s-}\!+\!Qv)\!+\!\Phi_s(Z_{s-})\Big]\tilde{N}(\dif v, \dif s)\bigg|^2\\
&=\mE\int_0^t\!\!\!\int_{|v|\leq1}\Big|\Phi_s^n(Z_s+Qv)-\Phi_s^n(Z_s)-\Phi_s\big(Z_s+Qv)+\Phi_s(Z_s)\Big|^2\nu(\dif v)\dif s\\
&\rightarrow0,\quad\text{as}\,\, n\rightarrow\infty.
\end{align*}
Finally, Krylov's estimate (\ref{kry2}) yields that
$$
\mE\left(\int_0^t\!\Big(B_n(s,Z_s)-B(s,Z_s)\Big)\dif s\right)\leq C\|B_n-B\|_{\mL^p_q(T)}\to 0, \quad \text{as}\,\,n\to\infty,
$$
which in turn as in $\cQ_4$ implies that
\begin{align*}
\cQ_1&\to \int_0^t\nabla_y \Phi_s(Z_s)\dif W_s,\quad\text{as}\,\, n\to\infty.
\end{align*}
Combing the above calculations, and noticing that $Z_s=\Phi^{-1}_s(\hat Z_s)$, we get the desired result.
\end{proof}

\subsection{Proof of main results}

Now we are in a position to give:

\begin{proof}[Proof of Theorem \ref{main1}]
By Lemma \ref{weak} and the classical Yamada-Watanabe principle \cite{Ya-Wa}, it is sufficient to show the pathwise uniqueness for solutions of
(\ref{sde1}). Put
\begin{align*}
\tau:=\inf\big\{t\geq 0:|L_t-L_{t-}|>1\big\}.
\end{align*}
Consider two solutions $Z_t^1$ and $Z_t^2$ of SDE (\ref{sde1}) both starting at $z\in\mR^{d_1+d_2}$. Then, by the interlacing property and the same argument
as \cite{Zh00}, we only need to show that
$$
Z_{t\wedge\tau}^1=Z_{t\wedge\tau}^2,\quad \text{a.s.}.
$$

Since uniqueness is a local concept, it is sufficient to focus on small $t\in[0,T]$. Let $\Phi_t(\cdot)$ given by (\ref{phi}) and define
$$
\hat Z^1_t:=\Phi_t(Z^1_t),\quad \hat Z^2_t:=\Phi_t(Z^2_t).
$$
According to Lemma \ref{zvon}, we have that $\hat Z_t^1-\hat Z_t^2$ satisfies the following SDE:
\begin{align*}
\hat Z_{t\wedge\tau}^1-\hat Z_{t\wedge\tau}^2
\!=\!\!&\int_0^{t\wedge\tau}\!\Big[\nabla_y \Phi_s\big(\Phi_s^{-1}(\hat Z^1_s)\big)-\nabla_y \Phi_s\big(\Phi_s^{-1}(\hat Z^2_s)\big)\Big]\dif W_s\\
&+\int_0^{t\wedge\tau}\Big[\tilde B(s,\hat Z^1_s)-\tilde B(s,\hat Z^2_s)\Big]\dif s\\
&+\!\!\int_0^{t\wedge\tau}\!\!\!\!\!\int_{|v|\leq 1}\!\!\Big[\Phi_s\big(\Phi_s^{-1}(\hat Z^1_s)\!+\!Qv\big)
\!-\!\hat Z^1_s\!-\!\Phi_s\big(\Phi_s^{-1}(\hat Z^2_s)\!+\!Qv\big)\!+\!\hat Z^2_s\Big]\tilde N(\dif v,\dif s).
\end{align*}
To shorten the notation, we denote by
\begin{align*}
\alpha(s)&:=\left|\nabla_y \Phi_s\big(\Phi_s^{-1}(\hat Z^1_s)\big)-\nabla_y \Phi_s\big(\Phi_s^{-1}(\hat Z^2_s)\big)\right|^2
+2\left\langle\hat Z_s^1-\hat Z_s^2,\tilde B(s,\hat Z^1_s)-\tilde B(s,\hat Z^2_s)\right\rangle,\\
\beta(s)&:=2\left\langle\hat Z_s^1-\hat Z_s^2,\nabla_y \Phi_s\big(\Phi_s^{-1}(\hat Z^1_s)\big)-\nabla_y \Phi_s\big(\Phi_s^{-1}(\hat Z^2_s)\right\rangle,\\
\xi(s,v)&:=\left|\Phi_s(\hat Z^1_s+Qv)-\Phi_s(\hat Z^1_s)-\Phi_s(\hat Z^2_s+Qv)+\Phi_s(\hat Z^2_s)\right|^2,\\
\eta(s,v)&:=\left|\hat Z_s^1-\hat Z_s^2\!+\!\Phi_s\big(\Phi_s^{-1}(\hat Z^1_s)\!+\!Qv\big)\!-\!\hat Z^1_s\!-\!\Phi_s\big(\Phi_s^{-1}(\hat Z^2_s)
\!+\!Qv\big)\!+\!\hat Z^2_s\right|^2\!-\!\big|\hat Z_s^1\!-\!\hat Z_s^2\big|^2.
\end{align*}
Now, It\^o's formula gives the following for any stopping time $\varsigma$,
\begin{align}
|\hat Z_{t\wedge\tau\wedge\varsigma}^1-\hat Z_{t\wedge\tau\wedge\varsigma}^2|^{2}&=\int_0^{t\wedge\tau\wedge\varsigma}\!\alpha(s)\dif s+\int_0^{t\wedge\tau\wedge\varsigma}\!\!\!\int_{|v|\leq 1}\xi(s,v)\nu(\dif v)\dif s\no\\
&\quad+\int_0^{t\wedge\tau\wedge\varsigma}\!\beta(s)\dif W_s+\int_0^{t\wedge\tau\wedge\varsigma}\!\!\!\int_{|v|\leq 1}\eta(s,v)\tilde N(\dif v,\dif s)\no\\
&= M_{t\wedge\tau\wedge\varsigma}+\int_0^{t\wedge\tau\wedge\varsigma}\!|\hat Z_s^1-\hat Z_s^2|^{2}\dif A_s, \label{00}
\end{align}
where $M_t$ is defined by
\begin{align*}
M_t&:=\int_0^{t}\beta(s)\dif W_s+\int_0^{t}\!\!\!\int_{|v|\leq 1}\eta(s,v)\tilde N(\dif v,\dif s),
\end{align*}
and $A_t$ is a continuous increasing process given by
\begin{align*}
A_t:=\int_0^t\frac{\alpha(s)}{|\hat Z_s^1-\hat Z_s^2|^2}\dif s+\int_0^{t}\!\!\!\int_{|v|\leq 1}\frac{\xi(s,v)}{|\hat Z_s^1-\hat Z_s^2|^2}\nu(\dif v)\dif s=:A_1(t)+A_2(t).
\end{align*}

It can be easily checked that $M_t$ is a martingale. Moreover, by (\ref{upd}), (\ref{es1}) and (\ref{w11}) we have that there is a constant $C_d>0$ such that
$$
\alpha(s)\leq C_d|\hat Z_s^1-\hat Z_s^2|^2\Big(1+\cM|\nabla^2_y U|(s,Z_s^1)+\cM|\nabla^2_y U|(s,Z_s^2)\Big),
$$
where $U$ satisfies (\ref{pde2}). Thus, we can deduce by (\ref{kry2}) and (\ref{mf}) that
$$
\mE A_1(t)\leq C+C\|\nabla^2_y U\|_{\mL^p_q(T)}<\infty.
$$

As for the term $A_2(t)$, let us denote by
$$
\cJ_v\Phi_s(z):=\Phi_s(z+Qv)-\Phi_s(z).
$$
Then, (\ref{es1}) yields that
\begin{align*}
\mE A_2(t)&\leq C\int_{|v|\leq 1}\!\mE\!\int_0^t\Big(\cM|\nabla\cJ_v\Phi_s|(Z^1_s)+\cM|\nabla\cJ_v\Phi_s|(Z^2_s)\Big)^2\dif s\nu(\dif v)\\
&\leq C\!\int_{|v|\leq 1}\!\|(\cM|\nabla\cJ_vU|)^2\|_{\mL^{p/2}_{q/2}(T)}\nu(\dif v)\\
&\leq C\!\int_{|v|\leq 1}\!\|\nabla\cJ_vU\|_{\mL^{p}_{q}(T)}^2\nu(\dif v).
\end{align*}
Our assumption $\beta>1-\tfrac{\alpha}{2}$ yields that
$$
2(\alpha/2+\beta)>2(\alpha+\beta-1)>\alpha.
$$
Consequently, it follows from \cite[Lemma 2.3]{Zh00} that
\begin{align*}
\mE A_2(t)\leq C\Big(\|U\|_{\mH^{\alpha+\beta,p}_q}^2+\|\nabla_yU\|_{\mH^{\alpha/2+\beta,p}_q}^2\Big)\!\int_{|v|\leq 1}|v|^{2(\alpha+\beta-1)}\nu(\dif v)<\infty.
\end{align*}
Therefore, $t\mapsto A(t)$ is a continuous strictly increasing process. Combing (\ref{00}), we get the desired result.
\end{proof}

With a little more efforts and as in \cite{Zh3}, we can give:

\begin{proof}[Proof of Corollary \ref{main11}]
For each $n\in\mN$, let $\chi_n^{d}(x)\in[0,1]$ be a nonnegative smooth function on $\mR^{d}$ with $\chi^{d}_n(x)=1$ for $x\in B_n:=\{x\in\mR^{d}: |x|\leq n\}$
and $\chi^d_n(x)=0$ for $x\notin  B_{n+1}$. Put
$$
F_n(t,x,y):=\chi^{d_1}_n(x)\chi^{d_2}_n(y)F(t,x,y),\quad G_n(t,x,y):=\chi^{d_1}_n(x)\chi^{d_2}_n(y)G(t,x,y).
$$
By Theorem \ref{main1}, for each $n\in\mN$, there exists a unique global strong solution $(X_t^n, Y_t^n)$ to SDE (\ref{sde1})
with coefficients $F_n, G_n$. For $n\geq k$, define
$$
\zeta_{n,k}:=\inf\{t\geq0: |X_t^n|+|Y_t^n|\geq k\}\wedge n.
$$
The uniqueness of the strong solution implies that
$$
\mP\Big((X_t^n,Y_t^n)=(X_t^k,Y_t^k),\,\forall t\in[0,\zeta_{n,k})\Big)=1,
$$
which yields that for $n\geq k$,
$$
\zeta_{k,k}\leq \zeta_{n,k}\leq \zeta_{n,n},\quad a.s..
$$
Put $\zeta_k:=\zeta_{k,k}$. We then have that $\{\zeta_k\}_{k\ge 1}$ is an increasing sequence of $(\sF_t)$-stopping times. And the following holds, for $n\geq k$,
$$
\mP\Big((X_t^n,Y_t^n)=(X_t^k,Y_t^k),\,\forall t\in[0,\zeta_{k})\Big)=1.
$$

Now, for each $k\in\mN$, we can define $X_t:=X^k_t$ and $Y_t:=Y_t^k$ for $t<\zeta_k$
and $\zeta:=\lim_{k\rightarrow\infty}\zeta_k$. It is easy to see that $(X_t, Y_t)$ is the unique strong solution to SDE (\ref{sde1})
up to the explosion time $\zeta$.
\end{proof}

\end{document}